\ifx\documentclass\undefined
\documentstyle[12pt]{article}
\else
\documentclass[12pt]{article}
\usepackage{latexsym}
\usepackage{amsfonts}
\usepackage{amsmath}
\usepackage{amssymb}
\usepackage{graphicx}
\usepackage{fullpage}
\fi

\author{ K\'aroly Bezdek\thanks{Partially supported by a Natural Sciences and 
Engineering Research Council of Canada Discovery Grant.}
and M\'arton Nasz\'odi\thanks{Partially supported by the Hung. Nat. Sci. Found.
(OTKA) grants K72537 and PD104744.}}

\font\tenBbb=msbm10 at 12pt         \font\sevenBbb=msbm9    \font\fiveBbb=msbm7
\newfam\Bbbfam
\textfont\Bbbfam=\tenBbb \scriptfont\Bbbfam=\sevenBbb
\scriptscriptfont\Bbbfam=\fiveBbb

\def\kkk{\null\hfill $\Box$\smallskip}

\newtheorem{theorem}{Theorem}[section]
\newtheorem{lemma}[theorem]{Lemma}

\newtheorem{problem}[theorem]{Problem}

\newtheorem{sled}[theorem]{Corollary}
\newtheorem{con}[theorem]{Conjecture}

\newtheorem{remark}[theorem]{Remark}

\newcommand{\proof}{{\noindent\bf Proof:{\ \ }}}

\newcommand{\Ee}{\mathbb E}
\newcommand{\Eh}{\Ee^3}
\renewcommand{\epsilon}{\varepsilon}

\newcommand{\B}{\mathbf B}
\newcommand{\BC}{{\B}(C)}
\newcommand{\VV}{\mathcal{V}}
\newcommand{\VT}{\VV^t}
\newcommand{\DD}{\mathcal{D}}
\newcommand{\DT}{\DD^t}

\newcommand{\NN}{\mathcal{N}}
\newcommand{\FF}{\mathcal{F}}
\newcommand{\GG}{\mathcal{G}}
\newcommand{\SC}{\mathcal{S}}
\newcommand{\st}{\colon}

\DeclareMathOperator{\inter}{int}
\DeclareMathOperator{\relint}{relint}
\DeclareMathOperator{\bd}{bd}

\DeclareMathOperator{\aff}{aff}

\newcommand{\conv}[1]{\left[ #1 \right]}

\title{Rigid Ball-Polyhedra in Euclidean $3$-Space 
}

\begin{document}
\maketitle
\date

\begin{abstract}
A ball-polyhedron is the intersection with non-empty interior of finitely many
(closed) unit balls in Euclidean $3$-space. One can represent the boundary of a
ball-polyhedron as the union of vertices, edges, and faces defined in a rather
natural way. A ball-polyhedron is called a simple ball-polyhedron if at every
vertex exactly three edges meet. Moreover, a ball-polyhedron is called a
standard ball-polyhedron if its vertex-edge-face structure is a lattice (with
respect to containment). To each edge of a ball-polyhedron one can assign an
inner dihedral angle and say that the given ball-polyhedron is locally rigid
with
respect to its inner dihedral angles if the vertex-edge-face structure of the
ball-polyhedron and its inner dihedral angles determine the ball-polyhedron up
to congruence locally. The main result of this paper is a Cauchy-type rigidity
theorem for ball-polyhedra stating that any simple and standard ball-polyhedron
is locally
rigid with respect to its inner dihedral angles.

\end{abstract}

\section{Introduction}

Let $\mathbb{E}^{3}$ denote the $3$-dimensional Euclidean space. As
in \cite{BN05} and \cite{BLNP07} a {\it ball-polyhedron} is the intersection
with non-empty interior of finitely many closed congruent balls in
$\mathbb{E}^{3} $. In fact, one may assume that the closed congruent
$3$-dimensional balls in question are of unit radius; that is, they are unit
balls of $\mathbb{E}^{3} $. Also, it is natural to assume that removing any of
the unit balls defining the intersection in question yields the intersection of
the remaining unit balls becoming a larger set. (Equivalently, using the
terminology introduced in \cite{BLNP07}, whenever we take a ball-polyhedron we
always assume that it is generated by a {\it reduced family} of unit balls.)
Furthermore, following \cite{BN05} and \cite{BLNP07} one can represent the
boundary of a ball-polyhedron in $\mathbb{E}^{3} $ as the union of {\it
vertices, edges}, and {\it faces} defined in a rather natural way. A {\it
standard ball-polyhedron} is one whose boundary structure is not
``pathological'', i.e., whose vertex-edge-face structure is a lattice (called the face lattice of the given standard ball-polyhedron). For the definitions of vertex, edge, face and standardness see Section~\ref{sec:combinatorics}. In order to get a more complete picture on ball-polyhedra, we refer the interested reader to \cite{BN05}, \cite{BLNP07} as well as \cite{KMP}.

One of the best known results in the geometry of convex polyhedra is Cauchy's rigidity theorem: If two convex polyhedra $P$ and $Q$ in $\mathbb{E}^{3}$ are combinatorially
equivalent with the corresponding faces being congruent, then the angles
between the corresponding pairs of adjacent faces are also equal and thus, $P$
is congruent to $Q$. Putting it somewhat differently the combinatorics of an
arbitrary convex polyhedron and its face angles completely determine its inner
dihedral angles. For more details on Cauchy's rigidity theorem and on its
extensions we refer the interested reader to \cite{Co93}. In our joint paper
\cite{BN05} we have been looking for analogues of Cauchy's rigidity theorem for
ball-polyhedra. In order to quote properly the relevant results from \cite{BN05}
we need to recall the following terminology. To each edge of a ball-polyhedron
in $\mathbb{E}^{3}$ we can assign an {\it inner dihedral angle}. Namely, take
any point $p$ in the
relative interior of the edge and take the two unit balls that contain
the two faces of the ball-polyhedron meeting along that edge. Now, the
inner dihedral angle along this edge is the angular measure of the intersection
of the two half-spaces supporting the two unit balls at $p$.
The angle in question is obviously independent of the choice of $p$.
Moreover, at each vertex of a face of a ball-polyhedron there is a
{\it face angle} formed by the two edges meeting at the given vertex
(which is, in fact, the angle between the two tangent halflines of the two edges
meeting at the given vertex). Finally, we say that the standard ball-polyhedron
$P$ in $\mathbb{E}^{3}$ is {\it globally rigid with respect
to its face angles} (resp., {\it its inner dihedral angles}) if the following
holds.
If $Q$ is another standard ball-polyhedron in $\mathbb{E}^{3}$ whose face
lattice is isomorphic
to that of $P$ and whose face angles (resp., inner dihedral angles)
are equal to the corresponding face angles (resp. inner dihedral angles)
of $P$, then $Q$ is congruent to $P$. 
We note that in \cite{BN05}, we used the word ``rigid'' for this notion. We
change that terminology to ``globally rigid'' because in the present paper we
consider a local version of the problem using the term ``locally rigid''.

Furthermore, a ball-polyhedron of
$\mathbb{E}^{3}$ is called {\it simplicial} if all its faces are bounded by
three edges. It is not hard to see that any simplicial ball-polyhedron is, in
fact, a standard one. Now, we are ready to state the main (rigidity) result of
\cite{BN05}: The face lattice and the face angles determine the inner dihedral
angles of any standard ball-polyhedron in $\mathbb{E}^{3}$. In particular, if
$P$ is a simplicial ball-polyhedron in $\mathbb{E}^{3}$, then $P$ is globally
rigid with respect to its face angles. The following fundamental analogue
question is still an open problem (see \cite{B10}, p. 63).

\begin{problem}\label{Bezdek}
Prove or disprove that the face lattice and the inner dihedral angles determine
the face angles of any standard ball-polyhedron in $\mathbb{E}^{3}$.
\end{problem}
 
One can regard this problem as an extension of the (still unresolved) conjecture
of Stoker \cite{St} according to which for convex polyhedra the face lattice and
the inner dihedral angles determine the face angles. For an overview on the
status of the Stoker conjecture and in particular, for the recent remarkable
result of Mazzeo and Montcouquiol on proving the infinitesimal version of the
Stoker conjecture see \cite{MM11}. The following special case of
Problem~\ref{Bezdek} has already been put forward as a conjecture in
\cite{BN05}. For this we need to recall that a ball-polyhedron is called a {\it
simple ball-polyhedron}, if at every vertex exactly three edges meet. Now, based
on our terminology introduced above the conjecture in question (\cite{BN05}, p.
257) can be phrased as follows.

\begin{con}\label{Bezdek-Naszodi}
Let $P$ be a simple and standard ball-polyhedron of $\mathbb{E}^{3}$. Then $P$
is globally rigid with respect to its inner dihedral angles.
\end{con}

We do not know whether the conditions of Conjecture~\ref{Bezdek-Naszodi} are
necessary. However, if the ball-polyhedron $Q$ fails to be a standard ball-polyhedron
because it possesses a pair of faces sharing more than one edge, then $Q$ is flexible (and so, it is not globally rigid) as shown in Section 4 of \cite{BN05}.

The main result of the present paper,
Theorem~\ref{Bezdek-Naszodi-rigidity-theorem}, is a local version of
Conjecture~\ref{Bezdek-Naszodi}.

\section{Main Result}

We say that the standard ball-polyhedron $P$ of $\mathbb{E}^{3}$ is {\it locally
rigid
with respect to its inner dihedral angles}, if there is an $\varepsilon>0$ with
the following property. If $Q$ is another standard ball-polyhedron of
$\mathbb{E}^{3}$ whose face lattice is isomorphic to that of $P$ and whose inner
dihedral angles are equal to the corresponding inner dihedral angles of $P$ such
that the corresponding faces of $P$ and $Q$ lie at Hausdorff distance at most
$\varepsilon$ from each other, then $P$ and $Q$ are congruent.

Now, we are ready to state the main result of this paper.

\begin{theorem}\label{Bezdek-Naszodi-rigidity-theorem}
Let $P$ be a simple and standard ball-polyhedron of $\mathbb{E}^{3}$. Then $P$
is locally rigid with respect to its inner dihedral angles.
\end{theorem}

Also, it is natural to say that the standard ball-polyhedron $P$ of
$\mathbb{E}^{3}$ is {\it locally rigid with respect to its face angles}, if
there is an
$\varepsilon>0$ with the following property. If $Q$ is another standard
ball-polyhedron of $\mathbb{E}^{3}$ whose face lattice is isomorphic to that of
$P$ and whose face angles are equal to the corresponding face angles of $P$ such
that the corresponding faces of $P$ and $Q$ lie at Hausdorff distance at most
$\varepsilon$ from each other, then $P$ and $Q$ are congruent. As according to
\cite{BN05} the face lattice and the face angles determine the inner dihedral
angles of any standard ball-polyhedron in $\mathbb{E}^{3}$ therefore
Theorem~\ref{Bezdek-Naszodi-rigidity-theorem} implies the following claim in a
straightforward way.

\begin{sled}
Let $P$ be a simple and standard ball-polyhedron of $\mathbb{E}^{3}$. Then $P$
is locally rigid with respect to its face angles.
\end{sled}

In the rest of this paper we give a proof of
Theorem~\ref{Bezdek-Naszodi-rigidity-theorem}.

\section{The Combinatorial Structure of a
Ball-Po\-ly\-hed\-ron}\label{sec:combinatorics}

Let $P$ be a ball-polyhedron in $\Eh$ given (as throughout the paper) by a
reduced family of generating balls.
A boundary point is called a {\it vertex} if
it belongs to at least three of the closed unit balls defining the
ball-polyhedron.
A {\it face} of the ball-polyhedron is the intersection of
one of the generating closed unit balls with the boundary of the
ball-polyhedron.
We say that the face of $P$ corresponds to the center of the generating ball.
Finally, if the intersection of two faces is non-empty, then it is the
union of (possibly degenerate) circular arcs. The non-degenerate
arcs are called {\it edges} of the ball-polyhedron.
Obviously, if a ball-polyhedron in $\mathbb{E}^{3}$ is generated by at least
three unit balls, then it possesses vertices, edges, and faces. Clearly, the
vertices, edges and faces of
a ball-polyhedron (including the empty set and the ball-polyhedron itself) are
partially ordered by inclusion forming the {\it vertex-edge-face structure} of
the given ball-polyhedron. 

We note that in \cite{BN05} the vertex-edge-face
structure of an arbitrary ball-polyhedron is incorrectly referred to as a face lattice.
Indeed, Figure~4.1 of \cite{BN05} shows an example of a ball-polyhedron whose
vertex-edge-face structure is not a lattice (with respect to inclusion). Thus, it is natural to define the following fundamental family of ball-polyhedra: a ball-polyhedron in $\mathbb{E}^{3}$ is a {\it standard ball-polyhedron} if its vertex-edge-face structure is a
lattice (with respect to inclusion). This is the case if, and only if, the
intersection of any two faces is either empty, or one vertex or one edge, and
every two edges share at most one vertex.
In this case, we simply call the vertex-edge-face structure in question the {\it
face lattice} of the standard
ball-polyhedron. This definition implies that any standard
ball-polyhedron of $\mathbb{E}^{3}$ is generated by at least four unit balls.

In connection with the above definition we note that the family of standard ball-polyhedra
was introduced and investigated in the more general, $n$-dimensional setting in \cite{BLNP07}. 
The $3$-dimensional case of that definition (Definition~6.4 in \cite{BLNP07}) coincides
with the definition given above. (See also Remark~9.1 and the paragraph preceding it
in \cite{BLNP07}.) For more insight on the vertex-edge-face structure of ball-polyhedra 
in $\mathbb{E}^{3}$ we refer the interested reader to \cite{KMP}.

\section{Infinitesimally Rigid Polyhedra, Dual Ball-Po\-ly\-hed\-ron, Trun\-ca\-ted
Delaunay Complex}\label{sec:infinitesimal}
\bigskip

In this section we introduce the notations and the main tools that are needed
for our proof of Theorem~\ref{Bezdek-Naszodi-rigidity-theorem}.

Recall that a {\it convex polyhedron} of $\mathbb{E}^{3}$ is a bounded
intersection of finitely many closed halfspaces in $\mathbb{E}^{3}$. A
\emph{polyhedral complex} in $\mathbb{E}^{3}$ is a finite family of convex
polyhedra such that 
any vertex, edge, and face of a member of the family is again a member of the
family, 
and the intersection of any two members is empty or a vertex or an edge or a
face of 
both members. In this paper a \emph{polyhedron} of $\mathbb{E}^{3}$ means the union of all members of a three-dimensional polyhedral complex in $\mathbb{E}^{3}$ possessing
the additional property that its (topological) boundary in $\mathbb{E}^{3}$ is a surface in $\mathbb{E}^{3}$ (i.e., a $2$-dimensional topological manifold embedded in $\mathbb{E}^{3}$).

We denote the convex hull of a set $C$ by $\conv{C}$. Following \cite{IS10}, we
call a polyhedron $Q$ in $\mathbb{E}^{3}$
\begin{itemize}
	\item \emph{weakly convex} if its vertices are in convex position (i.e.,
if its vertices are the vertices of a convex polyhedron);
	\item \emph{co-decomposable} if its complement in $\conv{Q}$ can be
triangulated (i.e., obtained as a simplicial complex) without adding new vertices;
	\item \emph{weakly co-decomposable} if it is contained in a convex
polyhedron $\tilde Q$ such that all vertices of $Q$ are vertices of $\tilde
Q$, 
	and the complement of $Q$ in $\tilde Q$ can be triangulated without
adding new vertices.
\end{itemize}

The boundary of every polyhedron in $\mathbb{E}^{3}$ is the disjoint union of
planar convex polygons and hence, it can be
triangulated without adding new vertices. Now, let $P$ be a polyhedron in
$\mathbb{E}^{3}$ and let $T$ be a triangulation of its boundary without adding
new vertices. We call the $1$-skeleton $G(T)$ of $T$ the {\it edge graph} of
$T$. By an {\it infinitesimal flex} of the edge graph $G(T)$ in $\mathbb{E}^{3}$
we mean an assignment of vectors to the vertices of $G(T)$ (i.e., to the
vertices of $P$) such that the displacements of the vertices in the assigned
directions induce a zero first-order change of the edge lengths: $(p_i-p_j)\cdot
(q_i-q_j)=0$ for every edge $p_ip_j$ of $G(T)$, where $q_i$ is the vector
assigned to the vertex $p_i$. An infinitesimal flex is called trivial if it is
the restriction of an infinitesimal rigid motion of $\mathbb{E}^{3}$. Finally,
we say that the polyhedron $P$ is {\it infinitesimally rigid} if every
infinitesimal flex of the edge graph $G(T)$ of $T$ is trivial. (It is not hard
to see that the infinitesimal rigidity of a polyhedron is a well-defined notion
i.e., independent of the triangulation $T$. For more details on this as well as
for an overview on the theory of rigidity we refer the interested reader to
\cite{Co93}.) We need the following remarkable rigidity theorem of Izmestiev and
Schlenker \cite{IS10} for the proof of
Theorem~\ref{Bezdek-Naszodi-rigidity-theorem}. 

\begin{theorem}\label{Izmestiev-Schlenker}(Izmestiev-Schlenker, \cite{IS10}) \\
Every weakly co-de\-com\-pos\-ab\-le polyhedron of $\mathbb{E}^{3}$ is infinitesimally rigid.
\end{theorem}

We note that Izmestiev and Schlenker \cite{IS10} give a different definition of
a polyhedron than ours, which yields a somewhat wider class of sets in $\Eh$.
Their theorem in its original form contains the additional restriction that  the
polyhedron is ``decomposable'' (i.e., it can be triangulated without new
vertices), which automatically holds for sets satisfying our narrower definition
of a polyhedron. Last but not least, one of the referees of our paper noted that by definition every weakly co-de\-com\-pos\-ab\-le polyhedron is in fact, a weakly convex one and therefore it is natural to state Theorem~\ref{Izmestiev-Schlenker} in the above form (i.e., not mentioning weakly convexity among the conditions).

The closed ball of radius $\rho$ centered at $p$ in $\mathbb{E}^{3}$ is denoted
by $\B(p,\rho)$. Also, it is convenient to use the notation $\B(p):=\B(p,1)$.
For a set $C\subseteq\mathbb{E}^{3}$ we denote the intersection of closed unit
balls with centers in $C$ by $\BC:=\cap\{\B(c)\st c\in C\}$. Recall that every
ball-polyhedron $P=\BC$ can be generated such that $\B(C\setminus\{c\})\neq \BC$
holds for any $c\in C$. Therefore whenever we take a ball-polyhedron $P=\BC$ we
always assume the above mentioned reduced property of $C$. The following duality
theorem has been proved in \cite{BN05} and it is also needed for our proof of
Theorem~\ref{Bezdek-Naszodi-rigidity-theorem}. 

\begin{theorem}\label{Bezdek-Naszodi-duality-theorem}(Bezdek-Nasz\'odi, \cite{BN05}) \\
Let $P$ be a standard ball-polyhedron of $\mathbb{E}^{3}$. Then the intersection
$P^*$ of the closed unit balls centered at the vertices of $P$ is another
standard ball-polyhedron whose face lattice is dual to that of $P$ (i.e., there
exists an order reversing bijection between the face lattices of $P$ and
$P^*$). 
\end{theorem}

For a more recent discussion on the above duality theorem and its generalizations we refer the interested reader to \cite{KMP}.



Let us give a detailed construction of the so-called {\it truncated
Delaunay complex} of an arbitrary ball-polyhedron, which is going to be the
underlying polyhedral complex of the given ball-polyhedron playing a central
role in our proof of Theorem~\ref{Bezdek-Naszodi-rigidity-theorem}.
We leave some of the proofs of the claims mentioned in the rest of this section
to the reader partly because they are straightforward and partly because they are also
well known (see \cite{AK00}, \cite{S91}, and in particular, \cite{EKS83}).

The \emph{farthest-point Voronoi tiling} corresponding to a finite set 
$C:=\{c_1,$ $\ldots,c_n\}$ in $\mathbb{E}^{3}$ 
is the family $\VV:=\{V_1,\ldots, V_n\}$ of closed {\it convex polyhedral sets}
$V_i:=\{x\in\mathbb{E}^{3}\st |x-c_i|\geq |x-c_j| \;\; {\rm for\ all}\  j\neq i,
1\le j\le n\}$, $1\le i\le n$. (Here a closed convex polyhedral set means a not
necessarily bounded intersection of finitely many closed halfspaces in
$\mathbb{E}^{3}$.)
We call the elements of $\VV$ \emph{farthest-point Voronoi cells}. 
In the sequel we omit the words ``farthest-point'' as we do not use the other
(more popular) Voronoi tiling: the one capturing closest points.

It is known that $\VV$ is a tiling of $\mathbb{E}^{3}$. We call the vertices,
(possibly unbounded) edges and (possibly unbounded) faces of the Voronoi cells
of $\VV$ simply the \emph{vertices}, \emph{edges} and \emph{faces} of $\VV$. 

The \emph{truncated Voronoi tiling} corresponding to $C$ is the family $\VT$ of
closed convex sets $\{V_1\cap\B(c_1),\ldots, V_n\cap \B(c_n)\}$. From
the definition it follows that
$\VT=\{V_1\cap P,\ldots, V_n\cap P\}$ where $P=\BC$.
We call elements of $\VT$ \emph{truncated Voronoi cells}.

Next, we define the (farthest-point) \emph{Delaunay complex} $\DD$ assigned to
the finite set $C=\{c_1,\ldots,c_n\}\subset\mathbb{E}^{3}$. 
It is a polyhedral complex on the vertex set $C$. For an index set
$I\subseteq\{1,\ldots,n\}$, the convex polyhedron $\conv{c_i\st i\in I}$ 
is a member of $\DD$ if, and only if, 
there is a point $p$ in $\cap\{V_i\st i\in I\}$ which is not contained in any
other Voronoi cell. 
In other words, $\conv{c_i\st i\in I}\in\DD$ if, and only if,
there is a point $p\in\mathbb{E}^{3}$ and a radius $\rho\geq 0$ such that 
$\{c_i\st i\in I\}\subset\bd \B(p,\rho)$ and $\{c_i\st i\notin I\}\subset\inter
\B(p,\rho)$. 
It is known that $\DD$ is a \emph{polyhedral complex}, in fact, it is a tiling
of $\conv{C}$ by convex polyhedra.

\begin{lemma}\label{lem:V}
	Let $C=\{c_1,\ldots,c_n\}\subset\mathbb{E}^{3}$ be a finite set, and
$\VV=\{V_1,\ldots,V_n\}$ be the corresponding Voronoi tiling of
$\mathbb{E}^{3}$. Then
	\begin{itemize}
		\item[(V)] 
				For any vertex $p$ of $\VV$, 
				there is an index set $I\subseteq\{1,\ldots,n\}$
with $\dim[c_i\st i\in I]=3$ 
				such that $\conv{c_i\st i\in I}\in\DD$ and
$p=\cap\{V_i\st i\in I\}$.\\
				And \emph{vice versa}: if
$I\subseteq\{1,\ldots,n\}$ with $\dim[c_i\st i\in I]=3$ is 
				such that $\conv{c_i\st i\in I}\in\DD$,  then
				$\cap\{V_i\st i\in I\}$ is a vertex of $\VV$.
		\item[(E)] 
				For any edge $\ell$ of $\VV$,
				there is an index set $I\subseteq\{1,\ldots,n\}$
with $\dim[c_i\st$ $i\in I]=2$ 
				such that $\conv{c_i\st i\in I}\in\DD$ and
$\ell=\cap\{V_i\st i\in I\}$.\\
				And \emph{vica versa}: if
$I\subseteq\{1,\ldots,n\}$ with $\dim[c_i\st i\in I]=2$ is 
				such that $\conv{c_i\st i\in I}\in\DD$,  then
				$\cap\{V_i\st i\in I\}$ is an edge of $\VV$.
		\item[(F)] 
				For any face $f$ of $\VV$,
				there is an index set $I\subseteq\{1,\ldots,n\}$
with $|I|=2$ 
				such that $\conv{c_i\st i\in I}\in\DD$ and
$f=\cap\{V_i\st i\in I\}$.\\
				And \emph{vica versa}: if
$I\subseteq\{1,\ldots,n\}$ with $|I|=2$ is 
				such that $\conv{c_i\st i\in I}\in\DD$, then
				$\cap\{V_i\st i\in I\}$ is a face of $\VV$.
	\end{itemize}
\end{lemma}

\proof
We outline the proof of (V) as the rest follows the same argument. Let $p$ be a
vertex of $\VV$, and let $I=\{i\in\{1,\ldots,n\}\st p\in V_i\}$. Now $p$ lies on
the boundary of some Voronoi cells. The centers corresponding to these Voronoi
cells are $\{c_i\st i\in I\}$. Since $p$ is shared by their Voronoi cells, these
centers are at an equal distance from $p$, in other words, they lie on a sphere
around $p$. Now, suppose that these centers are co-planar. Then, they lie on a
circle such that the line through the center of the circle, and perpendicular to
the plane of the circle, passes through $p$. Then all the Voronoi cells 
$\{V_i\st i\in I\}$ contain a relative neighborhood of $p$ within this line.
Thus, $p$ is not a vertex, a contradiction.


For the reverse statement:
Let $I$ be such that $\conv{c_i\st i\in I}\in\DD$ and $\dim[c_i\st i\in I]=3$.
It follows from the first condition on $I$ that $\cap\{V_i:i\in I\}\neq\emptyset$, 
and from that second condition that $\cap\{V_i:i\in I\}$ is a singleton, say $\{p\}$.
Clearly, $p$ is a vertex of $\VV$.
\kkk

\begin{figure}[hbp]
\unitlength1cm
\begin{minipage}[t]{.5\textwidth}
\begin{center}
\includegraphics[width=0.95\textwidth]{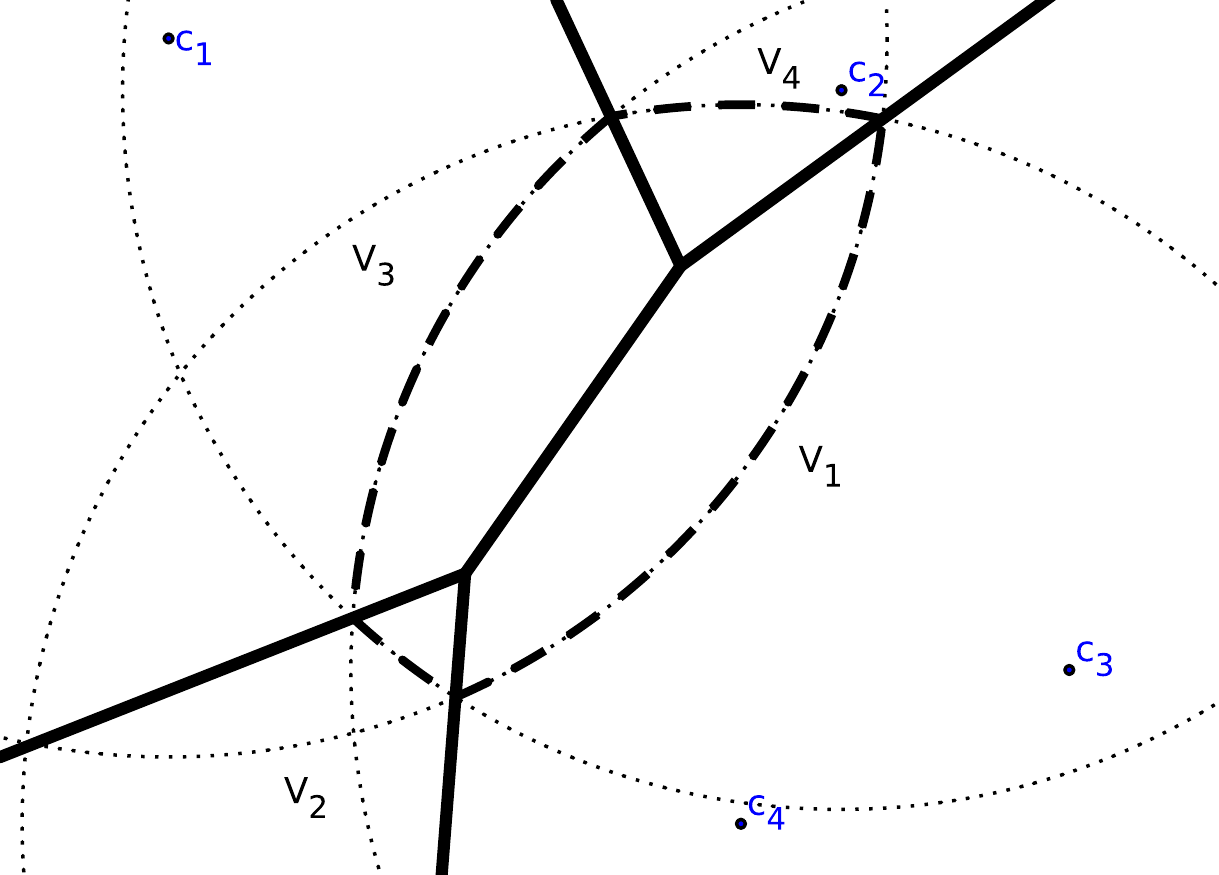}
\end{center}
\end{minipage}
\hfill
\begin{minipage}[t]{.5\textwidth}
\begin{center}
\includegraphics[width=0.95\textwidth]{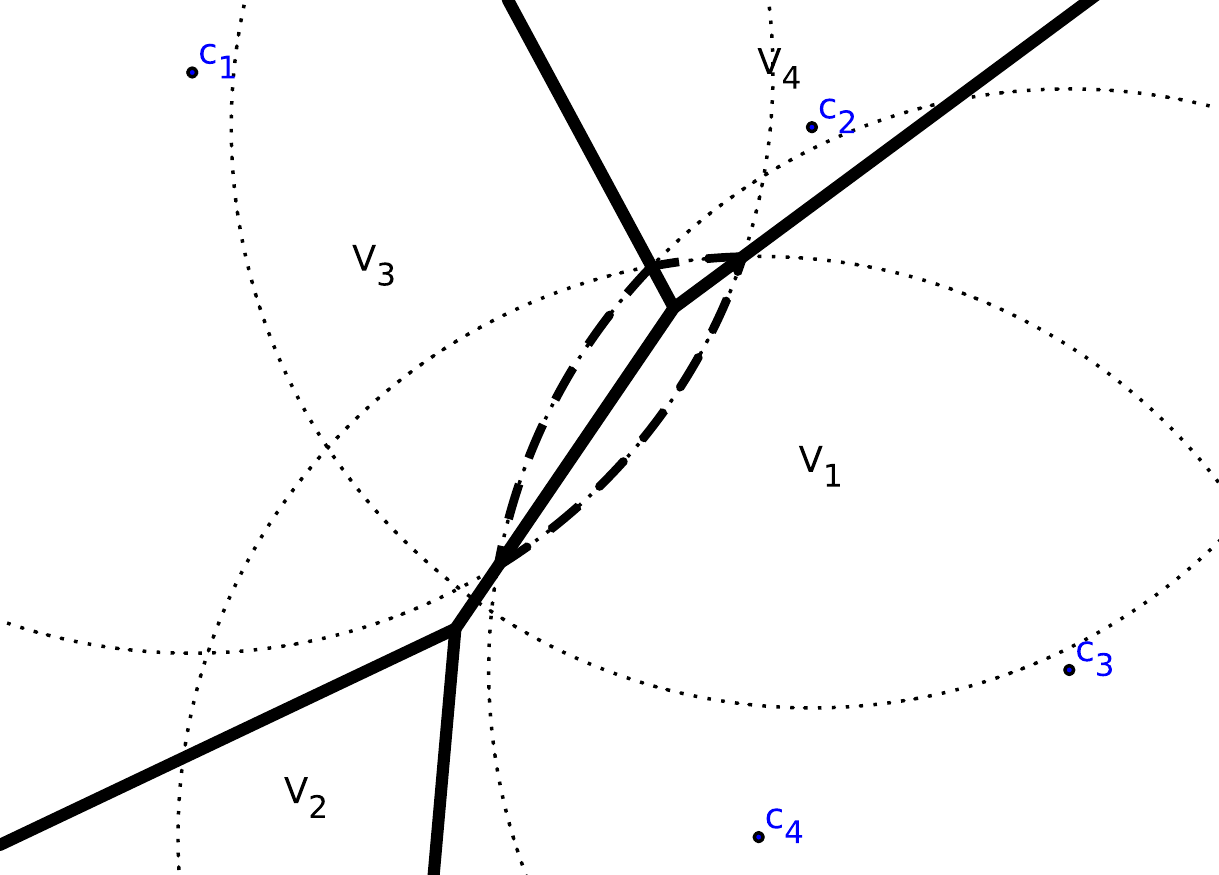}
\end{center}
\end{minipage}
\par\caption{
Given four points, $c_1,\ldots,c_4$. The bold solid lines bound the four Voronoi
cells,
$V_1,\ldots,V_4$. The bold dashed circular arcs bound the planar
ball-polyhedron -- a disk-polygon. The part of each Voronoi cell inside the
disk-polygon is the corresponding truncated Voronoi cell. 
On the first example,
$[c_1,c_3,c_4]$ and $[c_1,c_3,c_2]$ are the two-dimensional Delaunay cells, and
$[c_1,c_2], [c_1,c_3],[c_1,c_4],[c_2,c_3],[c_3,c_4]$ are the one-dimensional
Delaunay cells.
The truncated Delaunay complex coincides with the non-truncated one.
On the second example, the Voronoi and the Delaunay
complexes are the same as on the first, but the 
truncated Delaunay complex is different. The only two-dimensional truncated
Delaunay cell is $[c_1,c_3,c_4]$.
The one-dimensional truncated Delaunay cells are
$[c_1,c_3],[c_1,c_4],[c_3,c_4]$.
}\label{fig:voronoi1}
\end{figure}

We define the \emph{truncated Delaunay complex} $\DT$ corresponding to $C$
similarly to $\DD$:
For an index set $I\subseteq\{1,\ldots,n\}$, the convex polyhedron $\conv{c_i\st
i\in I}$ 
is a member of $\DT$ if, and only if, 
there is a point $p$ in $\cap\{V_i\cap\B(c_i)\st i\in I\}$ which is not
contained in any other truncated Voronoi cell. 
Note that the truncated Voronoi cells are contained in the ball-polyhedron
$\BC$.
Thus, $\conv{c_i\st i\in I}\in\DT$ if, and only if,
there is a point $p\in\BC$ and a radius $\rho\geq 0$ such that 
$\{c_i\st i\in I\}\subset\bd \B(p,\rho)$ and $\{c_i\st i\notin I\}\subset\inter
\B(p,\rho)$.

\section{Proof of Theorem~\ref{Bezdek-Naszodi-rigidity-theorem}}
\bigskip

\begin{lemma}\label{lem:novertex}
	Let $P=\BC$ be a simple ball-polyhedron in $\mathbb{E}^{3}$.
	Then no vertex of the Voronoi tiling $\VV$ corresponding to $C$ is on
$\bd P$, and no edge of $\VV$ is tangent to $P$.
\end{lemma}
\proof
By (V) of Lemma~\ref{lem:V}, at least four Voronoi cells meet in any vertex of
$\VV$. Moreover, the
intersection of each Voronoi cell with $\bd P$ is a face of $P$, since $P$ is
generated by a reduced set of centers. Hence, if a
vertex of $\VV$ were on $\bd P$ then at least four faces of $P$ would meet at a
point, contradicting the assumption that $P$ is simple.

Let $\ell$ be en edge of $\VV$, and assume that it contains a point $p\in\bd
P$. 
By the previous paragraph, $p\in\relint \ell$. From Lemma~\ref{lem:V}~(E) it
follows that $p$ is in the intersection of some Voronoi cells 
$\{V_i\st i\in I\}$ with $\dim[c_i\st i\in I]=2$. Clearly, $\ell$ is
orthogonal to the plane $\aff\{c_i\st i\in I\}$. Finally, there is an
$\varepsilon>0$ such that $P\cap\B(p,\varepsilon)=\B(\{c_i\st i\in
I\})\cap\B(p,\varepsilon)$ and hence,
$\ell$ must intersect $\inter P$, as $\ell$ intersects 
$\inter \big(\B(\{c_i\st i\in I\})\cap\B(p,\varepsilon)\big)$. \kkk

\begin{lemma}\label{lem:subcomplex}
	Let $P=\BC$ be a simple ball-polyhedron in $\mathbb{E}^{3}$.
	Then $\DT$ is a sub-polyhedral complex of $\DD$, that is
$\DT\subseteq\DD$, and faces, edges, and vertices of members of $\DT$ are again
members of $\DT$.
\end{lemma}
\proof
Clearly, $\DT\subseteq\DD$, and their vertex sets are identical (both are $C$).

First, we show that a (2-dimensional) face of a 3-dimensional member of $\DT$ is
again a member of $\DT$.
Let $\conv{c_i\st i\in I}\in\DT$ be a $3$-dimensional member of $\DT$. 
Then, the corresponding vertex (Lemma~\ref{lem:V}~(V)) $v$ of $\VV$ is in
$\inter P$ by Lemma~\ref{lem:novertex}.
For a given face of $\conv{c_i\st i\in I}$, there is a corresponding edge
(Lemma~\ref{lem:V}~(E)) $\ell$ of $\VV$. Clearly, $v$ is an endpoint of $\ell$.
Now, $\relint \ell\cap P\neq\emptyset$, and thus the face $\conv{c_i\st i\in I}$
of $\VV$ corresponding to $\ell$ is in $\DT$.

Next, let $\conv{c_i\st i\in I}\in\DT$ be a $2$-dimensional member of $\DT$ and
let $[c_i,c_j]$ be one of its edges. 
Then, for the corresponding edge $\ell$ of $\VV$ we have $\relint\ell\cap
P\neq\emptyset$. By Lemma~\ref{lem:novertex}, $\ell$ is not tangent to $P$, thus
$\relint\ell\cap \inter P\neq\emptyset$.
Now, $[c_i,c_j]$ corresponds to a face
(Lemma~\ref{lem:V}~(F)) $f$ of $\VV$. Clearly, $\ell$ is an edge of $f$.
Since an edge of $f$ intersects $\inter P$, we have $\relint f\cap P\neq\emptyset$
and hence, $f\cap P$ is a two-dimensional face of the truncated Voronoi tiling.
It follows that $[c_i,c_j]$ is in $\DT$. 
\kkk

The following lemma helps to understand the $2$-dimensional members of $\DT$.

Let $P=\BC$ be a simple and standard ball-polyhedron in $\mathbb{E}^{3}$.
Denote by $Q$ the polyhedral complex formed by the $3$-dimensional
members of $\DT$ and all of their faces, edges and vertices (i.e., we drop
``hanging'' faces/edges/vertices of $\DT$, that is, those faces/edges/vertices
that do not belong to a 3-dimensional member).
Clearly, $\cup Q$ is a subset of $\Eh$ and thus, its boundary is defined. We
equip this boundary with a polyhedral complex structure in the obvious way
as follows: we define the boundary of $Q$ as the collection of those faces,
edges and vertices of $Q$ that lie on the boundary of $\cup Q$. We denote this
polyhedral complex by $\bd Q$.

\begin{lemma}\label{lem:P}
Let $P=\BC$ be a simple and standard ball-polyhedron in $\mathbb{E}^{3}$, and
$Q$ be defined as above.
Then the $2$-dimensional members of $\bd Q$ are triangles, and a triangle
$\conv{c_1,c_2,c_3}$ is in $\bd Q$ if, and only if, the corresponding faces
$F_1,F_2,F_3$ of $P$ meet (at a vertex of $P$).
\end{lemma}

\proof
By Lemma~\ref{lem:subcomplex}, the $2$-dimensional members of $\bd Q$ are
$2$-dimensional members of $\DT$.
Let $\conv{c_i\st i\in I}\in\DT$ with $\dim[c_i\st i\in I]=2$. Then, clearly,
$\conv{c_i\st i\in I}\in\DD$ and, by Lemma~\ref{lem:V}~(E), it corresponds to an
edge $\ell$ of $\VV$ which intersects $P$. 
Now, $\ell$ is a closed line segment, or a closed ray, or a line. 
By Lemma~\ref{lem:novertex}, $\ell$ is not tangent to $P$, and (by
Lemma~\ref{lem:novertex}) $\ell$ has no endpoint on $\bd P$.
Thus, $\ell$ intersects the interior of $P$. We claim that $\ell$ has at least
one endpoint in $\inter P$. 
Suppose, it does not. Then $\ell\cap\bd P$ is a pair of points and so, the faces
of $P$ corresponding to indices 
in $I$ meet at more than one point. Since $|I|\geq 3$, it contradicts the
assumption that $P$ is standard.
We remark that this is a crucial point where we used the standardness of $P$.
So, $\ell$ has either one or two endpoints in $\inter P$. If it has two, then
the two distinct $3$-dimensional Delaunay cells 
corresponding to those endpoints (as in Lemma~\ref{lem:V}~(V)) are both members
of $\DT$ and
contain the planar convex polygon $\conv{c_I\st i\in I}$, and thus,
$\conv{c_I\st i\in I}$ is not on the boundary of $Q$.
If $\ell$ has one endpoint in $\inter P$, then there is a unique $3$-dimensional
polyhedron in $\DT$ 
(the one corresponding to that endpoint of $\ell$) that contains the planar
convex polygon $\conv{c_i\st i\in I}$. Moreover, in this case
$\ell$ intersects $\bd P$ at a vertex of $P$. Since $P$ is simple, that vertex
is contained in exactly three faces of $P$, and hence, 
$\conv{c_i\st i\in I}$ is a triangle. 

Next, working in the reverse direction, assume that $F_1,F_2$, and $F_3$ are
faces of $P$ that meet at a vertex $v$ of $P$. 
Then $v$ is in exactly three Voronoi cells, 
$V_1,V_2$ and $V_3$. Thus, $\conv{c_1,c_2,c_3}\in \DD$, and  $\ell:=V_1\cap
V_2\cap V_3$ is an edge of $\VV$. By the above argument,
$\ell$ has one endpoint in $P$ and so, $\conv{c_1,c_2,c_3}$ is a member of
$\DT$, 
and has the property that exactly one $3$-dimensional member of $\DT$ contains
it. It follows that $\conv{c_1,c_2,c_3}$ is in $\bd Q$.
\kkk

From the last paragraph of the proof and the fact that $P$ has at least one vertex,
we can deduce the following 
\begin{remark}
 With the notations and the assumptions of Lemma~\ref{lem:P}, $\DT$ contains
at least one $3$-dimensional cell, and the vertex set of $Q$ is $C$.
\end{remark}

We recall that the \emph{nerve} of a set family $\GG$ is the abstract simplicial
complex 
\[\NN(\GG):=\{\{G_i\in\GG\st i\in I\}\st \mathop\cap\limits_{i\in I}
G_i\neq\emptyset\}.\]

Now, let $P=\BC$ be a simple and standard ball-polyhedron in $\mathbb{E}^{3}$
and let $\FF$ denote the set of its faces. 
Let $\SC$ be the abstract simplicial complex on the vertex set $C$ generated by the 2-dimensional
members of $\bd Q$ (for the definition of $\bd Q$, see the paragraph preceding
Lemma~\ref{lem:P}), which are, according to Lemma~\ref{lem:P}, certain triples
of points in $C$. Both $\SC$ and the nerve $\NN(\FF)$ of $\FF$ are
$2$-dimensional abstract simplicial complexes. For the definition of an abstract
simplicial complex and its geometric realization, see \cite{Mat11}.

We claim that they both have the following ``edge property'': any edge is
contained in a $2$-dimensional simplex. Indeed, $\SC$
has this property by definition, since it is a simplicial complex generated by a
family of 2-dimensional simplices.
On the other hand, $\NN(\FF)$ also has this property, because $P$ is simple
and standard, and hence any edge of $P$ has a vertex as an endpoint which is a
point of intersection of three faces of $P$. 

Consider the mapping $\phi:c_i\mapsto F_i$ that maps each center point in $C$ to
the corresponding face of $P$. This is a bijection between the 0-dimensional
members of $\SC$ and the 0-dimensional members of $\NN(\FF)$.
By Lemma~\ref{lem:P} the 2-dimensional members of $\SC$ correspond via $\phi$ to
the 2-dimensional members of $\NN(\FF)$. By the ``edge property'' in the
previous paragraph, it follows that $\phi$ is an isomorphism of the two abstract
simplicial complexes, $\SC$ and $\NN(\FF)$.

By Theorem~\ref{Bezdek-Naszodi-duality-theorem}, $\NN(\FF)$ is
isomorphic to the face-lattice of another standard ball-po\-ly\-hed\-ron: $P^\ast$. 
Since $P^\ast$ is a convex body in $\mathbb{E}^{3}$ (i.e., a compact convex set
with non-empty interior in $\mathbb{E}^{3}$), the union of its faces is
homeomorphic to the $2$-sphere. Thus, $\SC$ as an abstract simplicial complex is
homeomorphic  to the $2$-sphere. On the other hand, $\bd Q$ is a geometric
realization of $\SC$. Thus, we have obtained that $\bd Q$ is a geometric
simplicial complex which is homeomorphic to the $2$-sphere. It follows that $Q$
is homeomorphic to the $3$-ball. So, we have that $Q$ is a polyhedron 
(the point being: it is topologically nice, that is, its boundary is a surface,
as required by the definition of a polyhedron in
Section~\ref{sec:infinitesimal}).

Clearly, $Q$ is a weakly convex polyhedron as $C$ is in convex position.
Furthermore, $Q$ is co-decomposable (and hence, weakly co-decomposable), as
$\DT$ is a sub-polyhedral complex of $\DD$ (by Lemma~\ref{lem:subcomplex}),
which is a family of convex polyhedra the union of which is $\conv{Q}=\conv{C}$.

So far, we have proved that $Q$ is a weakly convex and co-decomposable
polyhedron with triangular faces in $\mathbb{E}^{3}$. By
Theorem~\ref{Izmestiev-Schlenker}, $Q$ is infinitesimally rigid. Since $\bd Q$
itself is a geometric simplicial complex therefore its edge graph is rigid
because infinitesimal rigidity implies rigidity (for more details on that
see \cite{Co93}). Finally, we recall that the edges of the polyhedron $Q$ correspond
to the edges of the ball-polyhedron $P$, 
and the lengths of the edges of $Q$ determine (via a one-to-one mapping) 
the corresponding inner dihedral angles of $P$. It follows that $P$ is locally
rigid with respect to its inner dihedral angles.
\vskip0.5cm

{\bf Acknowledgements.} The authors wish to thank the anonymous referees for a number of helpful comments and suggestions.

\vspace{1cm}

\medskip

\noindent
K\'aroly Bezdek
\newline
Department of Mathematics and Statistics, University of Calgary, Canada,
\newline
Department of Mathematics, University of Pannonia, Veszpr\'em, Hungary,
\newline
Institute of Mathematics, E\"otv\"os University, Budapest, Hungary,
\newline
{\sf E-mail: bezdek@math.ucalgary.ca}

\medskip

\noindent
and

\medskip

\noindent
M\'arton Nasz\'odi
\newline
Institute of Mathematics, E\"otv\"os University, Budapest, Hungary,
\newline
{\sf E-mail: nmarci@math.elte.hu}

\end{document}